\documentclass[10pt]{article}
\usepackage{amsmath,amsfonts,amssymb,amsthm}
\usepackage{latexsym,xypic,eucal,mathrsfs,setspace}

\newcommand{\sa}{\emph{S.\ aureus}}
\newcommand{\ep}{\varepsilon}

\newcommand{\yone}{y_1}
\newcommand{\ytwo}{y_2}
\newcommand{\ythree}{y_3}

\newcommand{\kone}{k_1}
\newcommand{\ktwo}{k_2}

\newcommand{\xone}{x_1}
\newcommand{\xtwo}{x_2}

\title{Center Manifold and Lie Symmetry Calculations on a Quasi-chemical Model for Growth-death Kinetics in Food}

\author{Rachelle C. DeCoste    \\
             Department of Mathematical Sciences\\ The United States Military Academy\\ West Point, NY 10996 \\
              Tel.: 845-938-2530\\
              Rachelle.DeCoste@usma.edu   
        \and Louis  Piscitelle \\ U.S.\ Army RDECOM, Natick Soldier Center\\ Kansas Street, Natick, MA 01760\\
        Tel.: 508-233-4294\\
Louis.Piscitelle@us.army.mil  }

\newtheorem{theorem}{Theorem}[section]
\newtheorem{lemma}[theorem]{Lemma}

\begin{document}

\maketitle

\begin{abstract}
Food scientists at the U.S. Army's Natick Solider Center have developed a model for the lifecyle of the bacteria \emph{Staphylococcus aureus} in intermediate moisture bread.  In this article, we study this model using dynamical systems and Lie symmetry methods.  We calculate center manifolds and Lie symmetries for different cases of parameter values and compare our results to those of the food scientists.
\end{abstract}



\pagestyle{myheadings}
\thispagestyle{plain}
\markboth{R. DECOSTE AND L. PISCITELLE}{CENTER MANIFOLDS AND LIE SYMMETRIES}

\section{Introduction}

\subsection{The model}

To ensure the food that U.S. soldiers receive is as safe as possible, the growth of bacteria such as \emph{Staphylococcus aureus} (\sa) needs to be addressed.  The system of equations considered in this paper arises from a ``quasi-chemical'' kinetics model for the phases of the microbial life cycle of \sa \ in intermediate moisture bread.  The food scientists who developed the model confirmed its usefulness by fitting to it the data from observations on bread crumbs with varying conditions of water activity, pH and temperature.  The model differs from previous models in its attempt to model continuously the growth and death of the microorganism rather than focusing solely on either growth or inactivation.  The model was developed by food scientists at the Natick Soldier Center, see Taub, et al ~\cite{taub} for more on their techniques.

The model arose from the observation of four phases in the life cycle of \sa.  The cells pass through the various stages of metabolizing ($M$), multiplying ($M^*$), sensitization to death ($M^{**}$), and dead ($D$).  Additionally, the scientists hypothesized that there was an antagonist ($A$) present that would affect the cells.  They found that without this added element their original model did not fit the observed data with any accuracy.  The first step in the process describes cells moving from lag phase to growth phase ($M\rightarrow M^*$).  In the next step, cells multiply via binary division and then the newly multiplied cells interact with an antagonist ($M^*\rightarrow 2M^*+A$).  The last two steps represent two different pathways to death: the first with cells interacting with an antagonist, then passing to sensitization before death ($A+M^*\rightarrow M^{**}\rightarrow D)$ and lastly the cells experiencing natural death ($M^*\rightarrow D$).

The following equations represent the velocities of each of the above steps ($v$) as they relate to the concentrations of cells in various the phases.  Each equation has a rate constant ($k$) associated to it.
\begin{eqnarray}
v_1&=&k_1M\\
v_2&=&k_2M^*\\
v_3&=&(10^{-9})k_3M^*A\\
v_4&=&k_4M^*
\end{eqnarray}
Finally these velocities are represented by the following system of ordinary differential equations:
\begin{eqnarray}
\label{origfood1}\dot{M}&=&-v_1= -k_1M\\
\label{origfood2}\dot{M^*}&=&v_1+v_2-v_3-v_4=k_1M+M^*(G-\varepsilon A)\\
\label{origfood3}\dot{A}&=&v_2-v_3=M^*(k_2-\ep A)\\
\label{origfood4}\dot{D}&=&v_3+v_4 = M^*(k_4+\ep A)
\end{eqnarray}
where $G=k_2-k_4$ is the net natural growth rate and $\ep=10^{-9}k_3$.  It is assumed that all the rate constants have non-negative values.  The initial conditions at time zero are $M(0)=I$, the inoculum level $I\approx 10^3-10^4$,  and $M^*(0)=A(0)=D(0)=0$.

\subsection{A simplification}

We notice that the fourth equation is uncoupled since there are no terms involving the variable $D$ in any of the other equations and $\dot{D}$ depends on $M^*$ and $A$.  Therefore to investigate the dynamics of our system, we reduce to a system of three equations.  Renaming our variables ($\yone=M,\ \ytwo=M^*, \ \ythree=A$) we have the following system equivalent to equations \ref{origfood1}-\ref{origfood4}:
\begin{eqnarray}
\left(\begin{array}{c}\yone'\\ \ytwo' \\ \ythree' \end{array}\right)&=& \left(\begin{array}{ccc} -\kone & 0 & 0 \\ \kone & G & 0 \\ 0 & \ktwo & 0 \end{array}\right)\left(\begin{array}{c}\yone \\ \ytwo \\ \ythree \end{array}\right) + \left(\begin{array}{c}0 \\ -\ep\ytwo\ythree \\ -\ep\ytwo\ythree \end{array}\right)
\label{oureqn}
\end{eqnarray}

\section{Normalizing the system}
To consider the invariant manifold structure of a system, it is necessary to write the system in normal form as follows:
\begin{eqnarray}
\dot{x}&=&Ax+g(x,y)\\
\dot{y}&=&By+j(x,y)
\end{eqnarray}
with $(x,y)\in \mathbb{R}^n\times \mathbb{R}^m$,  the $n\times n$ matrix $A$ having eigenvalues with zero real part and the $m\times m$ matrix $B$ having eigenvalues with nonzero real part. The functions $g(x,y)$ and $j(x,y$) must be zero with zero first partial derivatives at the origin.  

The system (\ref{oureqn}) above is not in normal form since the $\ythree'$ equation corresponds to the zero eigenvalue piece and the nonlinear term of $\ythree'$ does not have all zero partial derivatives at the origin. Thus we must normalize by a change of coordinates using the eigenvectors of the matrix of the linear terms of the equation.  We will investigate the invariant manifolds in a neighborhood of $G=0$.  Writing our system in normal form for nonzero $G$ does not depend on the sign of $G$, so we treat the negative and positive case simultaneously.  We let $T$ be the matrix of eigenvectors of the eigenvalues of the linear terms of our system and let
\begin{eqnarray}
\left(\begin{array}{c}\yone \\ \ytwo \\ \ythree \end{array}\right)&=& T \left(\begin{array}{c} u \\ v \\ w \end{array}\right)=\left(\begin{array}{ccc} 0 & 0 & G+\kone \\ 0 & 1 & -\kone \\ 1 & \frac{\ktwo}{G} & \ktwo\end{array}\right)\left(\begin{array}{c} u \\ v \\ w \end{array}\right).
\end{eqnarray}
Using the inverse of the matrix $T$ we can solve for $u,\ v$ and $w$, find their derivatives and finally write our system in normal form as follows:
\begin{eqnarray}
u'&=&0\cdot u +\left(1-\frac{\ktwo}{G}\right)f(u,v,w)\\
\left(\begin{array}{c}v'\\ w'\end{array}\right)&=&\left(\begin{array}{cc} G & 0 \\ 0 & -\kone \end{array}\right) \left(\begin{array}{c}v\\w\end{array}\right) + \left(\begin{array}{c} f(u,v,w) \\ 0 \end{array}\right)
\end{eqnarray}
where $f(u,v,w)=-\frac{\ep}{G}(v-\kone w)(\ktwo v+G(u+\ktwo w))$.  Since $f(u,v,w)$ and its first partials with respect to $u, \ v,$ and $w$ are all zero at the origin $(u,v,w)=(0,0,0)$, we have our system in normal form and we see immediately that we have a one-dimensional center manifold in the case that $G\neq0$.  For $G>0$, we also have a one-dimensional stable and a one-dimensional unstable manifold.  For $G<0$, we have a two-dimensional stable manifold.  The system reduced to the center manifold simply becomes
\begin{equation}
u'=0.
\end{equation}

For the case $G=0$, we have a slightly simpler system of equations:
\begin{eqnarray}
\left(\begin{array}{c}\yone'\\ \ytwo' \\ \ythree' \end{array}\right)&=& \left(\begin{array}{ccc} -\kone & 0 & 0 \\ \kone & 0 & 0 \\ 0 & \ktwo & 0 \end{array}\right)\left(\begin{array}{c}\yone \\ \ytwo \\ \ythree \end{array}\right) + \left(\begin{array}{c}0 \\ -\ep\ytwo\ythree \\ -\ep\ytwo\ythree \end{array}\right)
\end{eqnarray}
Note that we now have two zero eigenvalues and one negative eigenvalue for the matrix in the linear term.  Since zero is a repeated eigenvalue, we must use generalized eigenvectors to find the normalization of this system.  Three such eigenvectors are $(0,0,1),\ (0,1,0),$ and $(1,-1,\frac{\ktwo}{\kone})$.  Then to transform our system we again let $T$ be the matrix consisting of these eigenvectors and let 
\begin{eqnarray}
\left(\begin{array}{c}\yone \\ \ytwo \\ \ythree \end{array}\right)&=& T \left(\begin{array}{c} u \\ v \\ w \end{array}\right)=\left(\begin{array}{ccc} 0 & 0 & 1 \\ 0 & 1 & -1 \\ 1 & 0 & \frac{\ktwo}{\kone}\end{array}\right)\left(\begin{array}{c} u \\ v \\ w \end{array}\right).
\end{eqnarray}
As above, this allows us to write our system in normal form:
\begin{eqnarray}
\left(\begin{array}{c}u'\\v'\end{array}\right)&=&\left(\begin{array}{cc}0 & \ktwo \\ 0 & 0 \end{array}\right)\left(\begin{array}{c}u\\v\end{array}\right)+\left(\begin{array}{c}g(u,v,w)\\g(u,v,w)\end{array}\right)\\
w'&=&-\kone w +0
\end{eqnarray}
where $g(u,v,w)=-\ep(v-w)(u+\frac{\ktwo}{\kone}w)$.  Since $g(u,v,w)$ is zero at the origin and all of its first partial derivatives are also zero at the origin, we can see that we have a two dimensional center manifold and a one dimensional stable manifold.
\section{Center manifold calculations}

Recall that a center manifold $W^c=\{(x,y)|y=h(x)\}$ is described by $h(x)$ 
where $h(0)=Dh(0)=0$.  We consider a system written in normal form
\begin{eqnarray}
\dot{x}&=&Ax+g(x,y)\\
\dot{y}&=&By+j(x,y)
\end{eqnarray}
with $A$ having eigenvalues with zero real part and $B$ eigenvalues with nonzero real part.  Then we determine $h(x)$ by finding the function that satisfies the following condition:
$$(Mh)(x)=Dh(x)[Ax+g(x,h(x))]-Bh(x)-j(x,h(x))=\overline{0}.$$
The sign of $G$ does not change the outcome of this calculation, thus we treat the case $G\neq 0$ at once.  We have $h:V\rightarrow \mathbb{R}^2$, $V\subset\mathbb{R}$ a neighborhood of the origin.  Thus let $h(x)=(h_1(x),h_2(x))=(ax^2+bx^3+O(x^4),cx^2+dx^3+O(x^4))$.  Then $f(x,h_1(x),h_2(x))=\ep(-a+\kone c)x^3+O(x^4)$ resulting in  
$$(Mh)(x)=\left(\begin{array}{c}-Gax^2+(-Gb-\ep a +\kone \ep c )x^3+O(x^4) \\ \kone c x^2+\kone d x^3 +O(x^4)\end{array}\right).$$
Solving for $(Mh)(x)=\overline{0}$, $h_1(x)=h_2(x)=O(x^4)$.  Thus up to third order, we have $h_1(x)=h_2(x)=0$, so a center manifold is simply the $u-$axis.

Next we consider the case $G=0$.  Here $h:V\rightarrow \mathbb{R}$, $V\subset\mathbb{R}^2$, a neighborhood of the origin. We let $h(x)=h(\xone,\xtwo)=a\xone^2+b\xtwo^2+c\xone\xtwo+d\xone^3+e\xtwo^3+f\xone^2\xtwo+j\xone\xtwo^2$.  Then we calculate
\begin{eqnarray*}
(Mh)(x)&=&\left(\begin{array}{cc}h_{\xone}(\xone,\xtwo), & h_{\xtwo}(\xone,\xtwo)\end{array}\right)\cdot\left(\begin{array}{c} \ktwo \xtwo + g(\xone,\xtwo,h(\xone,\xtwo)) \\ g(\xone,\xtwo,h(\xone,\xtwo)) \end{array}\right)+\kone h(\xone,\xtwo)\\
&=& (2a\ktwo+c\kone)\xone\xtwo+(a\kone)\xone^2+(c\ktwo+b\kone)\xtwo^2+(d\kone)\xone^3+(j\ktwo+e\kone)\xtwo^3\\
& &+(3d\ktwo-2a\ep-c\ep+f\kone)\xone^2\xtwo+(2f\ktwo-c\ep-2b\ep+j\kone)\xone\xtwo^2
\end{eqnarray*}
resulting in $h(\xone,\xtwo)=O(x^4)$, thus $h(\xone,\xtwo)=0$ up to order three.  Hence in this case the $uv-$plane is a center manifold.

\section{Lie Symmetry}

Recall that a Lie symmetry is a map from the set of solutions of a system of differential equations to the set itself.  For a system of first order ordinary differential equations 
\begin{equation} y_k'=\omega_k(t,\yone,\ytwo,\dots,y_n), \ \ k=1, \dots, n\label{difeq}\end{equation}
the Lie symmetries that transform the variables $t, y_1, \dots, y_n$ have infinitesimal generators of the form  \begin{equation} X=\xi\partial_t+\eta_1\partial_{y_1}+\eta_2\partial_{y_2}+\cdots+\eta_n\partial_{y_n}\end{equation}
where $\xi=\xi(t,y_1,y_2,\dots,y_n)$ and $\eta_k=\eta_k(t,y_1,y_2,\dots,y_n)$ for all $k$.  The infinitesimal generator must satisfy the Linearized Symmetry Condition:
\begin{equation}X^{(1)}(y_k'-\omega_k)=0, \ \ k=1,\dots,n \label{LSC}\end{equation}
when (\ref{difeq}) holds.  In this case the prolongation of $X$ is as follows:
\begin{equation} X^{(1)}=X+\eta_1^{(1)}\partial_{y_1'}+\eta_2^{(1)}\partial_{y_2'}+\cdots+\eta_n^{(1)}\partial_{y_n'}
\end{equation}
where $\eta_k^{(1)}$ is defined  as $\eta_k^{(1)}=D_t\eta_k-y_k'D_t\xi$.  The total derivative  $D_t$ in this case is $D_t=\partial_t+y_1'\partial_{y_1}+\cdots+y_n'\partial_{y_n}$.  Thus we have the following:
\begin{equation} \eta_{k}^{(1)}=\partial_t\eta_k+y_1'\partial_{y_1}\eta_k+y_2'\partial_{y_2}\eta_k+\cdots+y_n'\partial_{y_n}\eta_k-y_k'(\partial_t\xi+y_1'\partial_{y_1}\xi+y_2'\partial_{y_2}\xi+\cdots+y_n'\partial_{y_n}\xi). \end{equation}

A system of first order ODEs has an infinite number of symmetries.   We find symmetries by solving for the functions $\xi,\eta_k$ that satisfy the Linearized Symmetry Condition (\ref{LSC}).  This condition reduces to a system of PDEs which are computationally difficult to solve.  We use the ``Intro to Symmetry'' package in Mathematica and a script included in Cantwell \cite{cant} to calculate the symmetries for our system.  We are limited in the symmetries we can calculate by our computing power.    In the case $G\neq 0$ we calculate symmetries up to third order in our original coordinates $\yone,\ \ytwo,$ and $\ythree$ and then use a change of coordinates on our symmetries to rewrite in the coordinates $u,\ v,$ and $w$ of our equations in normal form.  Since the case $G=0$ involves simpler equations, we are able to calculate these symmetries directly from the equations in normal form, however we followed the same method as in the $G\neq 0$ case since we want to be able to compare cases.

\subsection{The case $G\neq 0$}
The infinitesimals of the Lie symmetries (up to order 3) are listed in an array with $\{\xi,\eta_1,\eta_2,\eta_3\}$, representing the infinitesimal generator $X=\xi\partial_t+\eta_1\partial_{\yone}+\eta_2\partial_{\ytwo}+\eta_3\partial_{\ythree}$.
\begin{eqnarray*}
X_1&=&\{1,0,0,0\}\\
X_2&=&\{\ytwo, -\kone\yone\ytwo, \kone\yone\ytwo+G\ytwo^2-\ep\ytwo^2\ythree,\ktwo\ytwo^2-\ep\ytwo^2\ythree\}\\
X_3&=&\{\ythree,-\kone\yone\ythree, \kone\yone\ythree+G\ytwo\ythree-\ep\ytwo\ythree^2,\ktwo\ytwo\ythree-\ep\ytwo\ythree^2\}\\
X_4&=&\{0,-\yone,\yone+\frac{G}{\kone}\ytwo-\frac{\ep}{\kone}\ytwo\ythree,\frac{\ktwo}{\kone}\ytwo-\frac{\ep}{\kone}\ytwo\ythree\}\\
X_5&=&\{\frac{1}{\kone}t,-t\yone,t\yone+\frac{G}{\kone}t\ytwo-\frac{\ep}{\kone}t\ytwo\ythree,\frac{\ktwo}{\kone}t\ytwo-\frac{\ep}{\kone}t\ytwo\ythree\}\\
X_6&=&\{\frac{-1}{\ep}\yone,\frac{\kone}{\ep}\yone^2,\frac{-\kone}{\ep}\yone^2-\frac{G}{\ep}\yone\ytwo+\yone\ytwo\ythree, \frac{-\ktwo}{\ep}\yone\ytwo+\yone\ytwo\ythree\}
\end{eqnarray*}
Then we transform the infinitesimal generators of the Lie symmetries found in the $y_i$ coordinates as follows.  If $X$ is an infinitesimal generator in $y_i$, then $\tilde{X}=(Xt)\partial_t+(Xu)\partial_u+(Xv)\partial_v+(Xw)\partial_w$ is the corresponding infinitesimal generator for a Lie symmetry in the $u,v,w$ coordinates 
The transformed symmetries in the form $\tilde{X}=\{\tilde{\xi},\tilde{\eta_1},\tilde{\eta_2},\tilde{\eta_3}\}$ where $\tilde{X}=\tilde{\xi}\partial_t+\tilde{\eta_1}\partial_u+\tilde{\eta_2}\partial_v+\tilde{\eta_3}\partial_w$: 
\begin{eqnarray*}
\tilde{X}_1&=&\{1,0,0,0\}\\
\tilde{X}_2&=&\{j(u,v,w),\frac{1}{G}(G-\ktwo)j(u,v,w)f(u,v,w),j(u,v,w)(Gv+f(u,v,w)),-\kone w j(u,v,w)\}\\
\tilde{X}_3&=&\{l(u,v,w), \frac{1}{G}(G-\ktwo)l(u,v,w)f(u,v,w),l(u,v,w)(Gv+f(u,v,w)),-\kone w l(u,v,w)\}\\
\tilde{X}_4&=&\{0,\frac{1}{G\kone}(G-\ktwo)f(u,v,w),\frac{1}{\kone}(Gv+f(u,v,w)),-w\}\\
\tilde{X}_5&=&\{\frac{t}{\kone},\frac{t}{G\kone}(G-\ktwo)f(u,v,w),\frac{t}{\kone}(Gv+f(u,v,w)),-tw\}\\
\tilde{X}_6&=&\{m(u,v,w),\frac{1}{G}(G-\ktwo)m(u,v,w) f(u,v,w),m(u,v,w)(Gv+f(u,v,w)), -\kone w m(u,v,w)
\end{eqnarray*}
where $f(u,v,w)$ is as above, $j(u,v,w)=v-\kone w$, $l(u,v,w)=u+\frac{\ktwo}{G}v+\ktwo w$ and $m(u,v,w)=-\frac{1}{\ep}(G+\kone)w$.

\subsection{The case $G=0$}
Again we calculate the infinitesimals of the Lie symmetries (up to order 3) of the original system with coordinates $\{\yone,\ytwo,\ythree\}$ and list them as $X=\{\xi,\eta_1,\eta_2,\eta_3\}$, representing the infinitesimal generator $X=\xi\partial_t+\eta_1\partial_{\yone}+\eta_2\partial_{\ytwo}+\eta_3\partial_{\ythree}$.
\begin{eqnarray*}
X_1&=&\{1,0,0,0\}\\
X_2&=&\{\ytwo,-\kone\yone\ytwo,\kone\yone\ytwo-\ep\ytwo^2\ythree,\ktwo\ytwo^2-\ep\ytwo^2\ythree\}\\
X_3&=&\{\ythree,-\kone\yone\ythree,\kone\yone\ythree-\ep\ytwo\ythree^2,\ktwo\ytwo\ythree-\ep\ytwo\ythree^2\}\\
X_4&=&\{0,\frac{\kone}{\ep}\yone,\frac{-\kone}{\ep}\yone+\ytwo\ythree,\frac{-\ktwo}{\ep}\ytwo+\ytwo\ythree\}\\
X_5&=&\{\frac{-1}{\ep}t,\frac{\kone}{\ep}t\yone, \frac{-\kone}{\ep}t\yone+t\ytwo\ythree,\frac{-\ktwo}{\ep}t\ytwo+t\ytwo\ythree\}\\
X_6&=&\{\frac{-1}{\ep}\yone,\frac{\kone}{\ep}\yone^2,\frac{-\kone}{\ep}\yone^2+\yone\ytwo\ythree,\frac{-\ktwo}{\ep}\yone\ytwo+\yone\ytwo\ythree\}
\end{eqnarray*}
Then we transform these to the $u,v,w$ coordinate system as above with $\tilde{X}=\{\tilde{\xi},\tilde{\eta_1},\tilde{\eta_2},\tilde{\eta_3}\}$ where $\tilde{X}=\tilde{\xi}\partial_t+\tilde{\eta_1}\partial_u+\tilde{\eta_2}\partial_v+\tilde{\eta_3}\partial_w$:
\begin{eqnarray*}
\tilde{X}_1&=&\{1,0,0,0\}\\
\tilde{X}_2&=&\{v-w,\frac{1}{\kone}(v-w)n(u,v,w),-\ep(v-w)p(u,v,w),-\kone(v-w)w\}\\
\tilde{X}_3&=&\{u+\frac{\ktwo}{\kone}w,\frac{1}{\kone}(u+\frac{\ktwo}{\kone}w)n(u,v,w),-\ep(u+\frac{\ktwo}{\kone}w)p(u,v,w),-\kone(u+\frac{\ktwo}{\kone}w)w\}\\
\tilde{X}_4&=&\{0, \frac{-1}{\ep\kone}n(u,v,w),p(u,v,w),\frac{\kone}{\ep}w\}\\
\tilde{X}_5&=&\{\frac{-t}{\ep},\frac{-1}{\ep\kone}t n(u,v,w),tp(u,v,w),\frac{\kone}{\ep}tw\}\\
\tilde{X}_6&=&\{\frac{-1}{\ep}w,\frac{-1}{\ep\kone}wn(u,v,w), wp(u,v,w),\frac{\kone}{\ep}w^2\}
\end{eqnarray*}
where $n(u,v,w)=\ep\ktwo w(-v+w)+\kone(\ktwo v+\ep u(-v+w))$ and $p(u,v,w)=(v-w)(u+\frac{\ktwo}{\kone}w)$.

\section{The connections between the center manifold and the Lie symmetry}
Recently Cicogna and Gaeta \cite{cico} have written about the connections between dynamical systems and Lie symmetries.  We are interested in particular in their results on invariant manifolds.  They have commented that any Lie symmetry of the system will leave invariant both the stable and unstable manifolds.  Due to the non-uniqueness of center manifolds, a Lie symmetry will map a center manifold to another (possibly the same) center manifold.  The following result indicates when a center manifold given by $\omega(u)$ will be invariant under a given Lie symmetry, in their notation $X=\phi\partial_u+\psi \partial_v$.

\begin{lemma}[Lemma 4 of \cite{cico} Chapter 7] A center manifold $w(u)$ is invariant under a Lie symmetry $X=\phi\partial_u+\psi \partial_v$ if and only if $$\psi(u,\omega(u))=(\partial_u(\omega(u)))\cdot \phi(u,\omega(u)).$$
\end{lemma}

For the case $G\neq 0$, $\omega(u)=\{0,0\}$ giving zero on the right side of this equality.  Thus the left side of this equation evaluated on the center manifold must always be zero if our center manifold is to be invariant under the action of the symmetry.  This is the case with all of our Lie symmetries as given above. 
For example consider $X_2$ with $\phi(u,v,w)=\frac{1}{G}(G-\ktwo)j(u,v,w)f(u,v,w)$ and $\psi(u,v,w)=\{j(u,v,w)(Gv+f(u,v,w)),-\kone w j(u,v,w)\}$.  Since $j(u,0,0)\equiv 0$,  $\psi(u,\omega(u))=\psi(u,0,0)=\{0,0\}$, thus satisfying the necessary and sufficient condition of the lemma.  It is easy to determine that the remainder of the symmetries in this case also leave the center manifold invariant.  Thus the center manifolds inherit these Lie symmetries.  However, in this case, since $v=w=0$, all of our symmetries become trivial.

Recall that in the case $G=0$ we found a center manifold to be the $uv-$plane.  Now, in the notation of our lemma, $\omega(u)=0$, and again the right side of our equation is zero.  Thus we must have $\psi(u,v,0)=0$ for any symmetry that leaves invariant this center manifold.  It can be easily checked to see that all of the symmetries listed above do indeed satisfy this necessary and sufficient condition.  In this case the center manifold again inherits the Lie symmetries which are now nontrivial.  The restriction of the system to our center manifold, the $uv-$plane, is 
\begin{eqnarray}
u'&=&\ktwo v - \ep u v \label{eqn51}\\
v'&=& -\ep u v \label{eqn52}.
\end{eqnarray}
The nontrivial symmetries inherited by this system are
\begin{eqnarray*}
\hat{X}_2&=&\{v,\frac{1}{\kone}vn(u,v,0),-\ep vp(u,v,0),0\}\\ 
\hat{X}_3&=&\{u,\frac{1}{\kone}un(u,v,0),-\ep up(u,v,w),0\}\\ 
\hat{X}_4&=&\{0, \frac{-1}{\ep\kone}n(u,v,0),p(u,v,0),0\}\\ 
\hat{X}_5&=&\{\frac{-1}{\ep}t,\frac{-1}{\ep\kone}t n(u,v,0),tp(u,v,0),0\}\\ 
\end{eqnarray*}

If we transform back to our original variables, we see that on the center manifold $u=\ythree$ and $v=\ytwo$, resulting in the system:
\begin{eqnarray}
\ytwo'&=&-\ep\ytwo\ythree\\
\ythree'&=&\ktwo\ytwo-\ep\ytwo\ythree
\end{eqnarray}
and the symmetries:
\begin{eqnarray}
\hat{X_2}&=&\ytwo\partial_t+\left(\ktwo\ytwo^2-\ep\ytwo^2\ythree\right)\partial_{\ytwo}-\ep\ytwo^2\ythree\partial_{\ythree}\\
\hat{X_3}&=&\ythree\partial_t+\left(\ktwo\ytwo\ythree-\ep\ytwo\ythree^2\right)\partial_{\ytwo}-\ep\ytwo\ythree^2\partial_{\ythree}\\
\hat{X_4}&=&\left(\frac{-\ktwo}{\ep}\ytwo+\ytwo\ythree\right)\partial_{\ytwo}+\ytwo\ythree\partial_{\ythree}\\
\hat{X_5}&=&\frac{-1}{\ep}t\partial_t+\left(\frac{-\ktwo}{\ep}t\ytwo+t\ytwo\ythree\right)\partial_{\ytwo}+t\ytwo\ythree\partial_{\ythree}
\end{eqnarray}

While we have calculated the infinitesimal generators, it would be interesting to determine the actual Lie symmetries on the center manifolds.  We would like to say precisely what these maps do to various trajectories on the center manifold and to the flow in general.  This is however, a very difficult question.  There is no known method that allows us to take the infinitesimal generators of any Lie symmetry and integrate them to find the actual symmetries.  The difficulty of this question is analogous to the solving of a system of  differential equations analytically. 

For example, if we consider $\hat{X}_4$ with $\eta_2(t,\ytwo,\ythree)=\frac{-\ktwo}{\ep}\ytwo+\ytwo\ythree$ and $\eta_3=\ytwo\ythree$, this means that, letting $\gamma$ be the parameter of the one-parameter Lie group, we need to solve the following for $\hat{\ytwo}$ and $\hat{\ythree}$, giving us the map $(\hat{\ytwo},\hat{\ythree})$ as our symmetry:
\begin{eqnarray}
\frac{d\hat{\ytwo}}{d\gamma}&=&\frac{-\ktwo}{\ep}\hat{\ytwo}+\hat{\ytwo}\hat{\ythree}\\
\frac{d\hat{\ythree}}{d\gamma}&=&\hat{\ytwo}\hat{\ythree}
\end{eqnarray}
This is equivalent to the system above.  Attempting to solve this system we find it equivalent to solving the following:
\begin{eqnarray}
\hat{\ytwo}&=& e^{\int (-\frac{\ktwo}{\ep}+\hat{\ythree})d\gamma}\\
\hat{\ythree}&=&e^{\int\hat{\ytwo}d\gamma}
\end{eqnarray}
with the initial conditions $\hat{\ytwo}(\gamma,\ytwo,\ythree)|_{\gamma=0}=\ytwo$ and $\hat{\ythree}(\gamma,\ytwo,\ythree)|_{\gamma=0}=\ythree$.

 This is something we continue to work on for this particular system as well as in general.

\subsection{Comparison to previous results}

Based on numerical solutions of the original system of equations Ross et.\ al \cite{ross} predicted trajectories for $M, \ M^*, A$ and $D$ with particular emphasis on the concentrations of $M^*$ (cells undergoing multiplication) and $A$ (the antagonist).  They found that the behavior depended on the values of the various constants $k_i$.  In particular, with $k_3=0$ and $G>0$, they found unrestrained growth of both $M^*$ and $A$.  For the values $k_3=0$ and a negative $G$, $A$ increases toward an upper limit and $M^*$ increases slightly but then begins to decrease toward zero.  For $k_3>0$ and $G>0$, both $M^*$ and $A$ increase initially, but then $M^*$ reaches a maximum and begins to decline while $A$ approaches an upper bound.  All of these analyses combined to indicate to the food scientists that the necessary constraints for growth-death kinetics are non-zero values for $k_3$ and positive values of $G$.

In our consideration of the system, we also found that $M^*(=y_3)$ and $A(=y_2)$ were the two variables that determined the behavior of the system.  In the $G\neq0$ case, the center manifold is the $u$-axis, which corresponds to $A$ when all other variables are zero, as on the center manifold.  
When $G = 0$, the reduced system on the center manifold is given by equations \ref{eqn51}  and \ref{eqn52}.  An inspection of this system, noting that $u = A$ and $v = M^*$, shows that the behavior is qualitatively identical to that found numerically in \cite{taub} for the case $k = [1\ 4\ 100\ 4]$, i.e. $k_2 = k_4 = 4$ resulting in $G = 0$.  In both the results are that $M^*$ goes to zero and $A$ approaches a constant value.

\section*{Acknowledgements}
This research was performed while the first author held a National Research Council Research Associateship Award jointly at the U.S. Army Natick Soldier Center, Natick, Massachusetts and the United States Military Academy, West Point, New York.


\begin{thebibliography}{10} 

\bibitem{cant}{\sc B.~J. Cantwell}, {\em Introduction to Symmetry Analysis}, Cambridge University Press, Cambridge, United Kingdom, 2002.

\bibitem{cico}{\sc G.~Cicogna and G.~Gaeta}, {\em Symmetry and Perturbation Theory in Nonlinear Dynamics}, Springer-Verlag, 1999.

\bibitem{ross} {\sc E.~Ross, I.~Taub, C.~Doona, F.~Feeherry, K.~Kustin}, {\em The mathematical properties of the quasi-chemical model for microorganism growth -- death kinetics in food}, International Journal of Food Microbiology, 99  (2005), pp. 157--171.

\bibitem{taub} {\sc I.~A. Taub, F.~E. Feeherry, E.~W. Ross, K.~Kustin, and C.~J.Doona}, {\em A Quasi-Chemical Kinetics Model for the Growth and Death of Staphylococcus aureus in Intermediate Moisture Bread}, Journal of Food Science, 68, No. 8 (2003), pp.~2530--2537.



\end{thebibliography}
\end{document}